\documentclass[a4paper]{iopart}
\usepackage{amssymb}
\usepackage{iopams}
\usepackage{setstack}
\usepackage{graphicx}
\usepackage{dcolumn}
\usepackage{bm}
\usepackage{hyperref}
\usepackage{color}
\usepackage{dsfont}
\usepackage{diagbox}
\usepackage{floatrow}
\usepackage[caption=false]{subfig}
\usepackage[export]{adjustbox}

\begin{document}

\title{Lower bounds for Ramsey numbers as a statistical physics problem}

\author{Jurriaan Wouters$^1$, Aris Giotis$^1$, Ross Kang$^2$, Dirk Schuricht$^1$ and Lars Fritz$^1$}

\address{$^1$ Institute for Theoretical Physics, Center for Extreme Matter and Emergent Phenomena, Utrecht University, Leuvenlaan 4, 3584 CE Utrecht, The Netherlands}
\address{$^2$ Department of Mathematics, Radboud University, Postbus 9010, 6500 GL Nijmegen, The Netherlands} \ead{d.schuricht@uu.nl, l.fritz@uu.nl}

%%%%%%%%%%%%%%%%%%%%%%%%%%%%%%%%%%%%
\begin{abstract}
Ramsey's theorem, concerning the guarantee of certain monochromatic patterns in large enough edge-coloured complete graphs, is a fundamental result in combinatorial mathematics.
In this work, we highlight the connection between this abstract setting and a statistical physics problem.
Specifically, we design a classical Hamiltonian that favours configurations in a way to establish lower bounds on Ramsey numbers. As a proof of principle we then use Monte Carlo methods to obtain such lower bounds, finding rough agreement with known literature values in a few cases we investigated. We discuss numerical limitations of our approach and indicate a path towards the treatment of larger graph sizes. 
\end{abstract}

\maketitle

%%%%%%%%%%%%%%%%%%%%%%%%%%%%%%%%%%%%
\section{Introduction}
%%%%%%%%%%%%%%%%%%%%%%%%%%%%%%%%%%%%
Ramsey numbers are a classic and fundamental optimisation parameter in the field of combinatorial mathematics~\cite{Ramsey30,graham1991}. To a lay audience this may be introduced via the so-called Party Problem: imagine that you are the host of a party, to which you invite some number of guests. However, you need to ensure, {\em no matter who attends}, that there are three guests that mutually know each other, or, alternatively, that there are three guests that are mutually stranger to one another. What is the least number of guests you must have for this condition to be met? The answer to this question is referred to as the Ramsey number $R(3,3)$ and it is of course easily reckoned: it is six. 

\begin{figure}[b]
\includegraphics[scale=0.8]{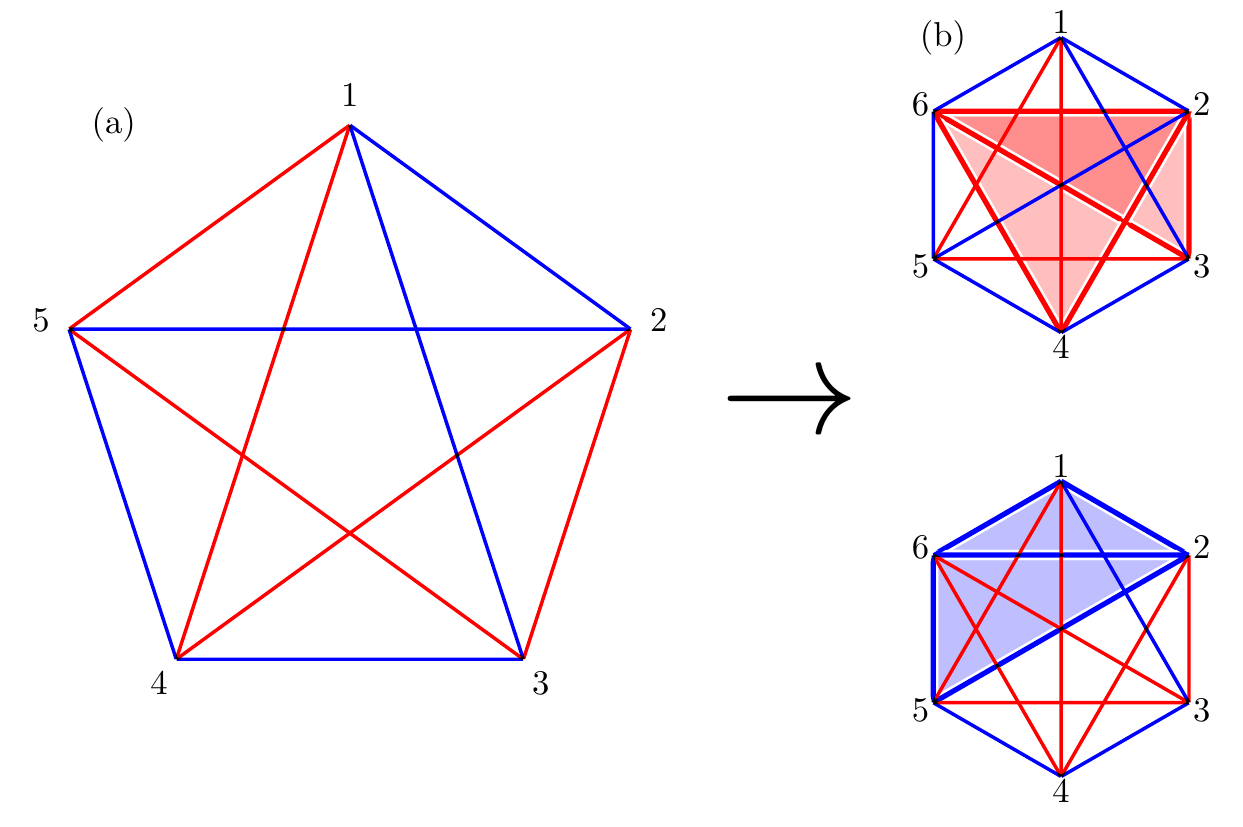}
\caption{(a) Complete graph with five vertices and a two-colour colouring without a monochromatic $3$-clique. (b) We add a vertex to the setup in (a). Irrespective of how we connect the additional vertex, we always end up with monochromatic $3$-cliques. Note that we only show two possible graphs.} \label{fig:partyproblem}
\end{figure}
More abstractly, one represents this by a complete graph, to each of its edges of which are individually assigned exactly one of two colours, red or blue. (A complete graph is a graph in which there is an edge between every possible pair of vertices.) See  Figure~\ref{fig:partyproblem} for an illustration. For the Party Problem, we interpret each vertex as a guest and one of the colours, blue, say, as signifying that the corresponding pair of guests are already acquainted, while the other colour signifies that the pair are unacquainted. In Figure~\ref{fig:partyproblem}(a) we represent one specific party with five guests. Three guests mutually knowing each other corresponds to a $3$-clique or {\em triangle} all of whose edges are blue, whereas three guests being mutually stranger occurs as a red triangle. (A $k$-clique is a complete $k$-vertex subgraph.) For the conditions of the Party Problem, we can check for each triangle if its edges are all blue or all red, that is, if it forms a {\em monochromatic} $3$-clique. Observe in Figure~\ref{fig:partyproblem}(a) that there is no monochromatic $3$-clique, which certifies that having five guests is not enough, i.e.~that the Ramsey number $R(3,3)$ is strictly greater than five. In Figure~\ref{fig:partyproblem}(b), we increase the number of vertices to six, where two possible ways to colour the resulting graph are shown. It happens that both of them contain monochromatic $3$-cliques. In fact, by a careful organisation, one can check by hand all the cases to prove that every red/blue-edge-colourings of a complete $6$-vertex graph possesses a monochromatic $3$-clique, i.e.~confirming that the answer to the Party Problem is indeed six~\cite{greenwood_combinatorial_1955}. Alternatively, one can `brute force' enumerate and check with a computer all $2^{15}$ possible red/blue-colourings. 

As we alluded to above, for lower bounds to certain Ramsey numbers, it suffices to find some edge-colouring configuration on a complete graph that is absent of monochromatic cliques of some given size. Now, viewing the vertices as lattice sites and the coloured edges as `interactions' between suitable degrees of freedom, the task becomes reminiscent of a statistical physics setup. Finding what we in short refer to as a `clique-free configuration' is then tantamount to identifying a certain state. It is the aim of this article to make this line of reasoning more precise. To this end, we first discuss a suitably chosen statistical physics system for which finding the ground state becomes equivalent to finding a clique-free configuration. We then apply Monte Carlo techniques~\cite{NewmanBarkema99}  to actually identify such clique-free configurations, in turn providing lower bounds for the corresponding Ramsey number. We note that a similar approach was followed by Gaitan and Clark~\cite{GaitanClark12} and mentioned by Exoo~\cite{Exoo1989OnTC}, who related Ramsey numbers to a combinatorial optimisation problem. This was then treated using an adiabatic quantum evolution algorithm, both in simulations~\cite{GaitanClark12} and experimentally implemented using Squid qubits~\cite{Bian-13}. An overview of other numerical efforts for Ramsey numbers is given by Haanp\"a\"a~\cite{haanpaa_computational_2007}. We contend that the reminder of this interesting link could spur further developments and progress for Ramsey numbers, especially in light of work on frustrated spin systems.

This article is organised as follows. In the next section we give a more formal definition of the Ramsey numbers. In Section~\ref{sec:main} we elaborate on our line of reasoning and set up the statistical physics system. This is then treated using Monte Carlo methods as discussed in Section~\ref{sec:Monte}, the results of which are presented in Section~\ref{sec:results}. The challenges to the Monte Carlo approach are discussed in Section~\ref{sec:challange}. Details of our numerical implementation are given in the appendix.

%%%%%%%%%%%%%%%%%%%%%%%%%%%%%%%%%%%%%%%%%%%%%%%%%%%%%%%%%%%%%%%
\section{Ramsey number}
%%%%%%%%%%%%%%%%%%%%%%%%%%%%%%%%%%%%%%%%%%%%%%%%%%%%%%%%%%%%%%%
We consider a more general {\em multicolour} Ramsey number, $R(x_1,\ldots,x_l)$. It is the smallest positive integer $k$ for which, on a complete graph of size $k$ each of whose edges is coloured with one of $l$ different colours, there is guaranteed to be a clique of size $x_i$ that is monochromatic in colour $i$, for at least one of the values $i=1,\dots,l$. As mentioned, the Party Problem corresponds to $R(3,3)=6$. Ramsey~\cite{Ramsey30} introduced and proved that these parameters are all well-defined.

Even though the Ramsey problem is easy to state and grasp, the precise value of $R(x_1,\ldots,x_l)$ is unknown except for the smallest choices of parameters. Over the decades, many bounds, asymptotic limits and relations between Ramsey numbers have been pursued and obtained; see the survey of Radziszowski~\cite{Radziszowski21} for an updated overview. For instance, while $R(4,4)=18$ is known exactly, even $R(5,5)$ is only known to lie in the range from 43 to 48\footnote{There is a famous anecdote about Erd\H{o}s~\cite{Spencer94}. He imagined a superior alien race landing on Earth and threatening to destroy it if humankind fails to deliver the value for $R(5,5)$. In that case, he would advise to redirect all computational resources and mathematicians to this task. If, instead, they asked for $R(6,6)$, he would advise rather to focus on destroying the aliens.}. Selected results relevant to our work are given in (the brackets in) Tables~\ref{tab:2color} and~\ref{tab:multicolor}. 

A na\"ive strategy to gather information about Ramsey numbers is through brute force search algorithms. If in this way one finds a clique-free configuration on $N$ vertices, it implies that the corresponding Ramsey number $R$ satisfies $N<R$. The obvious problem with this is that the number of possible configurations grows tremendously with $N$, making the search infeasible even for relatively modest values of $N$. 
With a more guided search, one might arrive more rapidly at such desired configurations. Our hope is that modern methods in statistical physics could help to devise such searches, and open a path towards better methods for finding lower bounds.

%%%%%%%%%%%%%%%%%%%%%%%%%%%%%%%%%%%%%%%%%%%%%%%%%%%%%%%%%%%%%%%
\section{Linking to a statistical physics problem}\label{sec:main}
%%%%%%%%%%%%%%%%%%%%%%%%%%%%%%%%%%%%%%%%%%%%%%%%%%%%%%%%%%%%%%%
The key step is to translate the problem of finding clique-free configurations into the task of identifying ground states of a suitably chosen statistical physics system. To this end, we define an energy functional that assigns an energy to each possible colouring of a complete graph. The functional is constructed such that monochromatic cliques are energetically unfavourable; see below for our specific choice. Then a suitable technique, in our case the Monte Carlo Metropolis algorithm, is used to minimise the energy, ideally finding a clique-free configuration. While this does not precisely determine the value of the Ramsey number, it can certify a lower bound.  
 
The search will be guided towards clique-free configurations, if the employed energy functional satisfies the following conditions: (i) assign zero energy to clique-free configurations and (ii) be strictly positive for all configurations containing monochromatic cliques (of the given sizes). Obviously, these conditions can be adapted, for example, by shifting the overall energy. 

To be specific, let us consider the Ramsey number $R(x_1,\ldots,x_l)$ on a complete graph with $N$ vertices. For each colour $i$ we write all $x_i$-cliques as $Q_j^{x_i}$, $j=1,\ldots, n_{x_i}$, where $n_{x_i}={N \choose x_i}$ is the total number of such cliques. To each $x_i$-clique $Q_j^{x_i}$ we then assign an energy via
\begin{equation}\label{eq:localenergy}
\epsilon_i(Q^{x_i}_j)=\left\{\begin{array}{ll} 1& \mathrm{if\;all\;edges\;in\;}Q_j^{x_i}\mathrm{\;have\;colour\;}i,\\ 0& \mathrm{otherwise}.\end{array}\right.
\end{equation}
The energy of a configuration on the whole graph is then given by summing over all colours $i$ as well as all possible cliques: 
\begin{equation}\label{eq:energy}
E=\sum_{i=1}^l\left[K_i\sum_{j=1}^{n_{x_i}}\epsilon_i(Q^{x_i}_j)\right], \quad K_i>0.
\end{equation}
Thus every monochromatic clique of colour $i$ contributes $K_i$ to the energy. The requirement $K_i>0$ ensures that condition (ii) is met. Empirically we determined that $K_i\propto \frac1{x_i}$ results in the fastest convergence. A similar energy functional (cost function) was used in the quantum algorithm approach~\cite{GaitanClark12,Bian-13}. Of course, the choice~(\ref{eq:energy}) is by no means unique, thus opening a path for future improvements of this approach. 

%%%%%%%%%%%%%%%%%%%%%%%%%%%%%%%%%%%%%%%%%%%%%%%%%%%%%%%%%%%%%%%
\section{Method}\label{sec:Monte}
%%%%%%%%%%%%%%%%%%%%%%%%%%%%%%%%%%%%%%%%%%%%%%%%%%%%%%%%%%%%%%%
Having defined an energy functional and thus a statistical physics model, we can use the methods of statistical mechanics to study its properties. The number of configurations (microstates) grows exponentially with the number of vertices $N$, for example as $\sim 2^{N(N-1)/2}$ in the two-colour case. This implies that already for 12 vertices we have to deal with about $7.4 \times10^{19}$ different microstates\footnote{Permutations reduce the number of inequivalent microstates compared to $2^{N(N-1)/2}$ because of indistinguishable vertices.}. A brute force search algorithm, consequently, is relying on luck.  For illustration, performing a brute force search for 12 vertices, we obtained a configuration without (4,4)-cliques on average after about 60 seconds. For comparison, the Monte Carlo algorithm (see Section~\ref{sec:MCA} below) discussed below required mere milliseconds. 

%%%%%%%%%%%%%%%%%%%%%%%%%%%%%%%%%%%%%%%%%%%%%%%%%%%%%%%%%%%%%%%
\subsection{Monte Carlo simulations}\label{sec:MCA}
%%%%%%%%%%%%%%%%%%%%%%%%%%%%%%%%%%%%%%%%%%%%%%%%%%%%%%%%%%%%%%%
Monte Carlo methods~\cite{NewmanBarkema99} are one of the most popular techniques to study statistical physics problems. These  methods allow one to tackle large system sizes by identifying a representative set of microstates. 

In our case we wish to find configurations that minimise the energy functional in~(\ref{eq:energy}). We do so by employing the well-known Monte Carlo Metropolis algorithm in combination with simulated annealing. We start with a random configuration of coloured edges and calculate its energy  $E_{\rm{initial}}$. Next, we randomly pick a edge, change its colour (blue to red, red to blue, blue to green, etc.), and calculate the energy $E_{\rm{final}}$ of the resulting configuration. If $E_{\rm{final}}<E_{\rm{initial}}$, we accept the new configuration, i.e.~we continue the algorithm with it. If $E_{\rm{final}}>E_{\rm{initial}}$, we accept the new configuration with probability $\exp[(E_{\rm{initial}}-E_{\rm{final}})/T]$, where $T>0$ is the temperature. Then we proceed by changing the colour of the next edge and iterate this procedure over all edges. Since we are more interested in swapping edges of existing monochromatic cliques, our algorithm visits these edges with higher frequency; for details, see~\ref{app:biasedbonds}. Usually in simulated annealing, the temperature $T$ is lowered gradually, to variationally approach the thermodynamic ground state. 

Due to the long-range interactions encoded in the energy functional~(\ref{eq:energy}), the energy landscape turns out to be rather erratic, i.e.~there are transitions between configurations with large energy differences by changing the colour of a single edge. This is different from the situation usually encountered when studying physical systems. It implies that the energy landscape shows many local minima, with the algorithm showing a tendency to get stuck in. Therefore, on top of the annealing procedure, we vary the temperature to allow the system to escape these local minima. For details, see~\ref{app:temperature}. We note that our ultimate goal is to find a zero-energy state (i.e.~clique-free configuration) for a given number of vertices, not calculating thermodynamic quantities, so these non-physical alterations are permitted. There is, however, one word of caution in order. The energy functional possesses a large number of degenerate configurations in all energy sectors, including the ground state manifold. This implies that the model shares some similarities with frustrated spin systems but also spin glass problems, which limits the applicability of the Metropolis algorithm. We leave it to the future to further investigate the relation between the present model and these fields of research. 

%%%%%%%%%%%%%%%%%%%%%%%%%%%%%%%%%%%%%%%%%%%%%%%%%%%%%%%%%%%%%%%
\subsection{Adjacency matrix}
%%%%%%%%%%%%%%%%%%%%%%%%%%%%%%%%%%%%%%%%%%%%%%%%%%%%%%%%%%%%%%%
At the implementation level, we store the information about the edges in the so-called adjacency matrix $A$. For its definition, we choose an ordering $1,\ldots,N$ of the vertices like in Figure~\ref{fig:partyproblem}. Then, $A$ is an $N\times N$-matrix, where the  matrix element $A_{pq}$ corresponds to the edge between vertex $p$ and vertex $q$. Obviously, the adjacency matrix is symmetric, $A_{pq}=A_{qp}$. In the case of the two-colour Ramsey number we choose the convention that a blue edge corresponds to $A_{pq}=+1$, whereas a red edge is encoded as $A_{pq}=-1$. The diagonals are not used during the algorithm; we set them to $A_{pp}=0$. To illustrate this, the configuration shown in Figure~\ref{fig:partyproblem}(a) for $R(x_1,x_{-1})=R(3,3)$ is represented by the $5 \times 5$-matrix
\begin{eqnarray}
A=\left( \begin{array}{ccccc} 0 & 1& 1 &-1 &-1 \\ 1&0&-1&-1&1\\ 1&-1&0&1&-1 \\ -1&-1&1&0&1\\ -1&1&-1&1&0   \end{array}\right).
\end{eqnarray}
The adjacency matrix gives us access to the colouring of a $3$-clique (triangle) $Q^3_j=\{p,q,r\}$ via its elements $A_{pq}$, $A_{qr}$, and $A_{pr}$. To determine whether the triangle is monochromatic, we have to compare these elements, i.e.~if $A_{pq}=A_{qr}=A_{pr}=\pm 1$ we have $\epsilon_{\pm 1}\left(Q^3_j\right)=1$ and $\epsilon_{\mp 1}\left(Q^3_j\right)=0$, whereas we have $\epsilon_{\pm 1}\left(Q^3_j\right)=0$ in all other cases.

It is straightforward to generalise the adjacency matrix to $l$-coloured graphs for $R(x_1,\ldots,x_l)$. The off-diagonal elements, representing the coloured edges, take values $1,\ldots,l$. For larger cliques $x_i>3$, the condition in~(\ref{eq:localenergy}) generalises to $A_{pq}=i$ for all $p,q\in\{p_1,\ldots,p_{x_i}\}$ belonging to the $x_i$-clique $Q_j^{x_i}$. In this notation, a Metropolis step thus involves switching an edge variable, meaning $A_{pq}$ changes its colour. We have implemented several optimisations in the determination of the clique energy which are listed in~\ref{app:forloopreduction} and turn out to reduce the computation time significantly.

%%%%%%%%%%%%%%%%%%%%%%%%%%%%%%%%%%%%%%%%%%%%%%%%%%%%%%%%%%%%%%%
\subsection{Recycled initial state}
%%%%%%%%%%%%%%%%%%%%%%%%%%%%%%%%%%%%%%%%%%%%%%%%%%%%%%%%%%%%%%%
There is a further way to reduce the computation time. Suppose the algorithm has identified a clique-free configuration with $N$ vertices, thus already showing that $R(x_1,\ldots,x_l)\ge N+1$. Now we use this clique-free configuration when adding another vertex. Typically, the resulting configuration on the $(N+1)$-vertex graph will contain significantly fewer monochromatic cliques than a randomly generated state (since the initial $N$-clique was clique-free), thus the Monte Carlo procedure will usually converge much faster.

\begin{figure}[t]
	\centering
	\subfloat[\label{fig:R_5_5_10}]{\includegraphics[scale=0.5,valign=b]{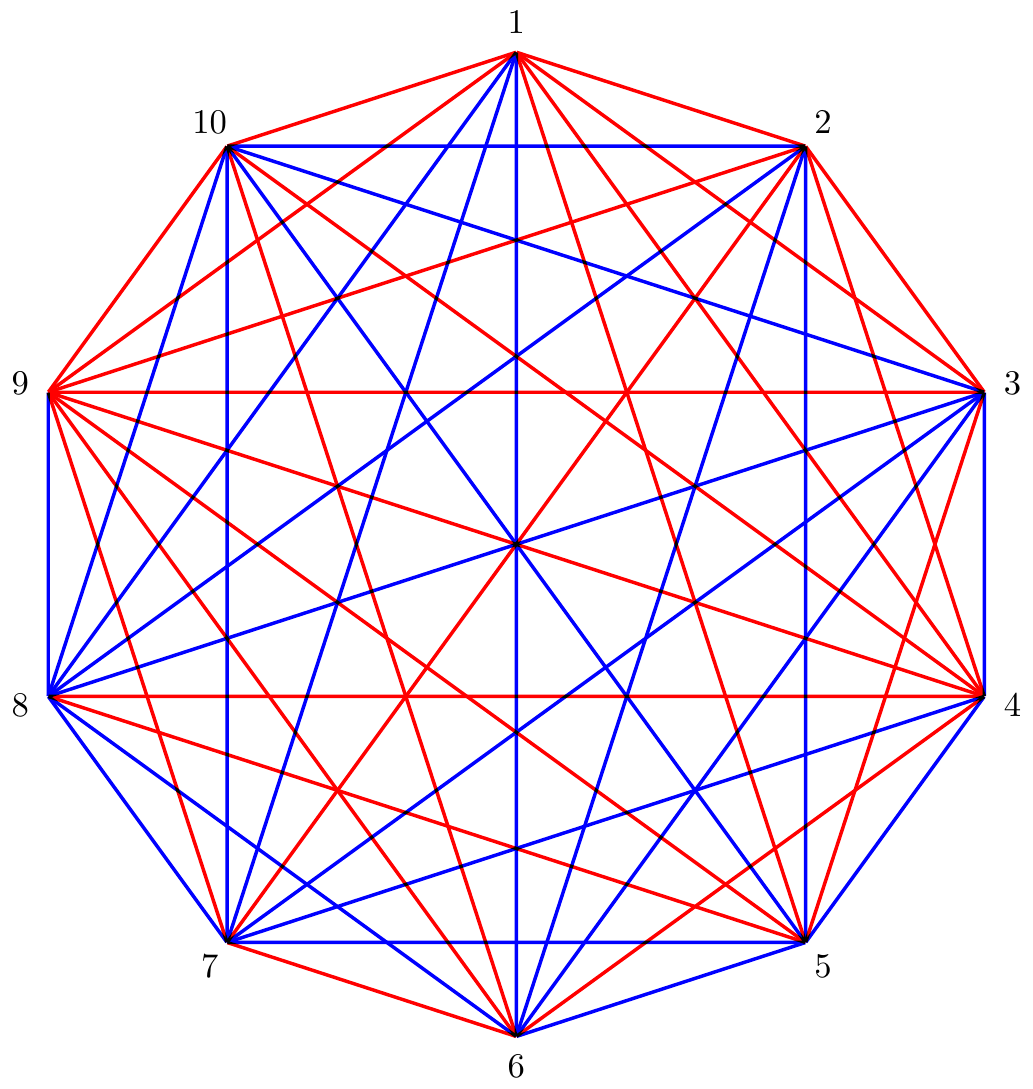}}
	\hfill
	\subfloat[\label{fig:R_5_5_10_with_extension}]{\includegraphics[scale=0.5,valign=b]{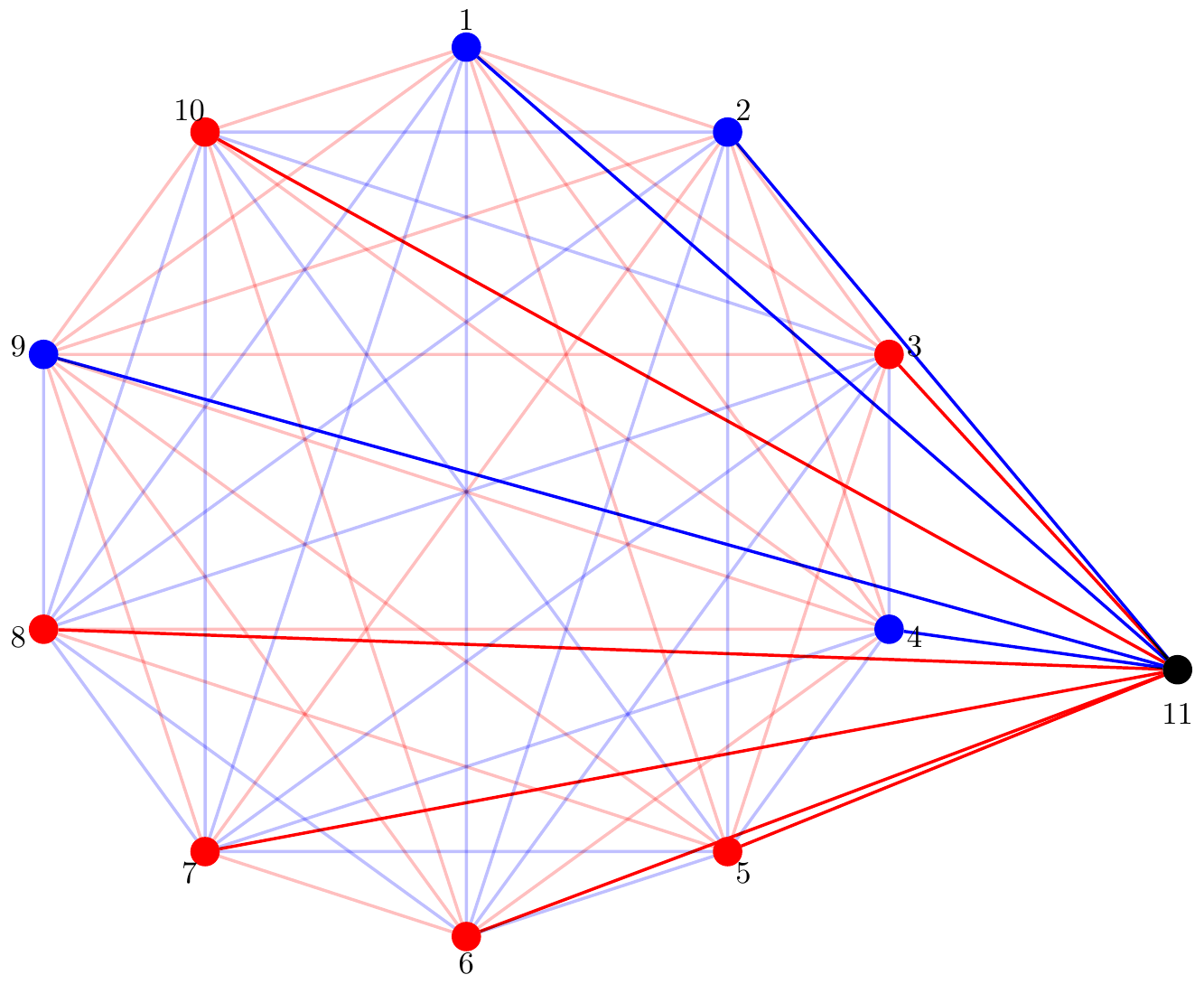}}
	\caption{An example of a 10-vertex non-monochromatic 5-clique configuration (a). Adding an additional vertex in (b), we choose the initial colours of the emerging edges such that the resulting initial configuration is less likely to contain monochromatic cliques (see main text). Note that (b) is generally not clique-free.} \label{fig:R55}
\end{figure}

Furthermore, there is a simple pragmatic way to choose the extra bonds for the new $N+1$-vertex state, which again increases the chances of getting an initial state with few monochromatic cliques. Consider for example the configuration for $R(5,5)$ on a 10-vertex graph shown in Figure~\ref{fig:R55}(a). For each of the vertices we can read off the least occurring colour (from the adjacency matrix), visualised by the colour of the vertices in Figure~\ref{fig:R55}(b). Now, when connecting the 'old' vertices to the additional one (vertex 11), we choose the respective edges to have the least dominant colour (see the colouring of the additional edges). This reduces the chance of monochromatic cliques in the initial configuration.

%%%%%%%%%%%%%%%%%%%%%%%%%%%%%%%%%%%%%%%%%%%%%%%%%%%%%%%%%%%%%%%
\section{Results}\label{sec:results}
%%%%%%%%%%%%%%%%%%%%%%%%%%%%%%%%%%%%%%%%%%%%%%%%%%%%%%%%%%%%%%%
After discussing our implementation of the Monte Carlo simulation, we now turn to the obtained results. For a given number of vertices $N$ the algorithm aims at finding a clique-free configuration, which in turn provides a strict lower bound for the corresponding Ramsey number, $N<R$. 

\begin{table}[t]
	\begin{tabular}{r|c|c|c|c|c|c}
		\diagbox{$x_1$}{$x_2$}	& 3		& 4			& 5				& 6 			& 7 		& 8\\
		\hline
		3 				& 6 (6)	& 9 (9)		& 14 (14) 		& 18 (18) 		& 23 (23)	& 28 (28)\\
		\hline
		4 				&		& 18 (18)	& 25 (25) 		& 34 (36-41)	& &\\
		\hline
		5 				&		& 			& 41 (4$3$-48)	&&&\\
		
	\end{tabular}
	\caption{Lower bounds obtained for the two-colour Ramsey numbers $R(x_1,x_2)$ using the Monte Carlo simulations to obtain clique-free configurations. In brackets we state the known values (or respective ranges)~\cite{Radziszowski21}.}\label{tab:2color}
\end{table}
As a proof of principle we used this approach to consider the most accessible Ramsey numbers, i.e.~ones for which the lower bound is (likely) to be below 60 vertices. The lower bounds obtained from our simulations are presented in Tables~\ref{tab:2color} and~\ref{tab:multicolor}. For the simplest cases, e.g.~$R(4,4)=18$, exact values are known. For larger cliques and most multicoloured graphs only the ranges are known, e.g.~$45\le R(3,3,5)\le 57$. For comparison, we list the known literature values~\cite{Radziszowski21} in the brackets.
\begin{table}[t]
	\begin{tabular}{c | c | c |c | c}
$R(3,3,3)$	& $R(3,3,4)$	& $R(3,3,5)$ 	& $R(3,4,4)$	& $R(3,3,3,3)$\\
		\hline
17 (17)		& 29 (30) 		& 43 (45-57)	& 48 (55-77)	& 43 (51-62)	
	\end{tabular}
\caption{Lower bounds obtained for the multi-colour Ramsey numbers $R(x_1,\ldots,x_l)$ using the Monte Carlo simulations to obtain clique-free configurations. In brackets we state the known values (or respective ranges)~\cite{Radziszowski21}.}\label{tab:multicolor}
\end{table}

\begin{figure}[t]
	\centering
	\includegraphics[width=0.8\textwidth]{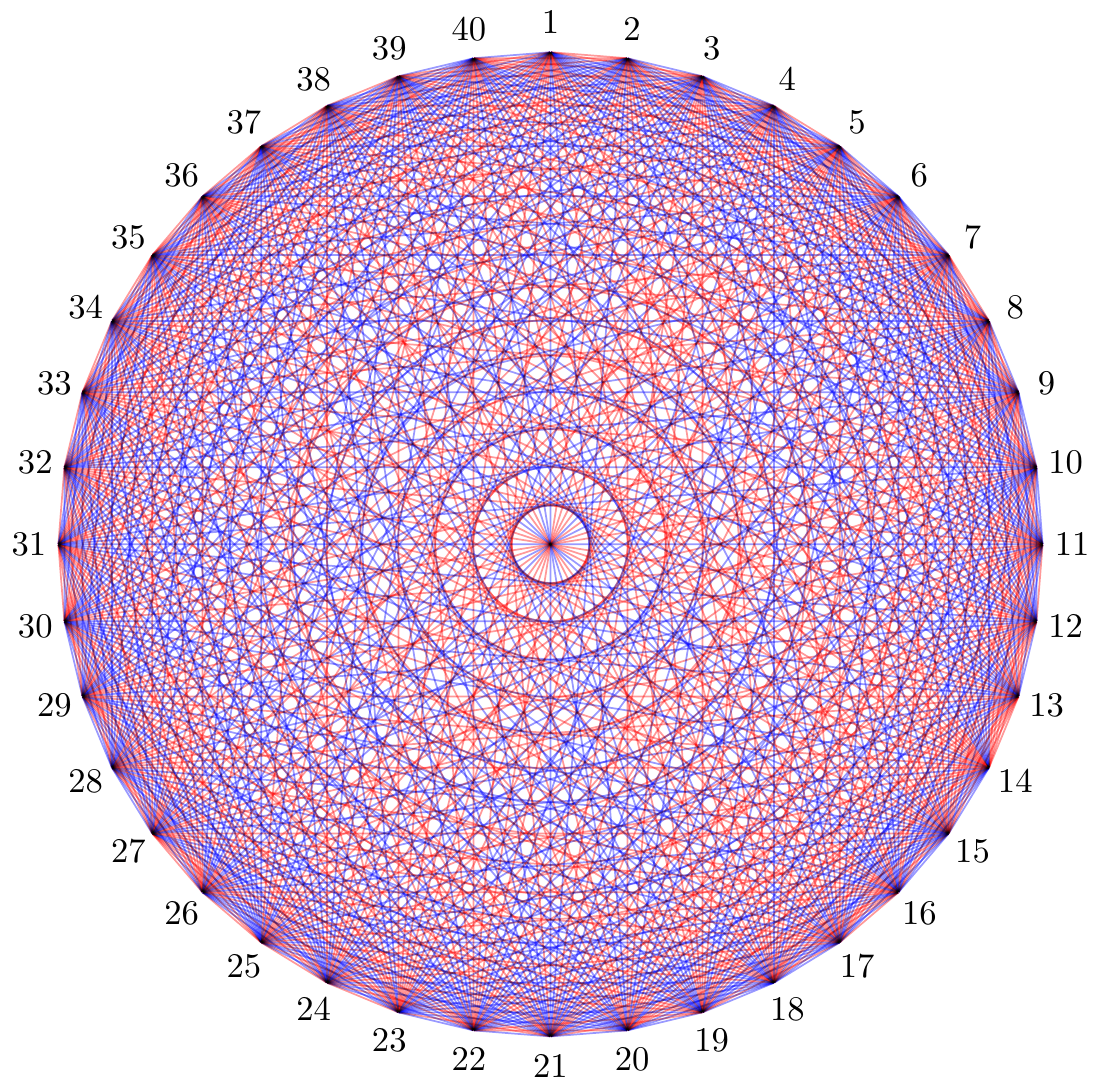}
	\caption{clique-free configuration on a complete two-colour graph with $N=40$ vertices. The existence of this configuration implies the strict lower bound $40<R(5,5)$ for the corresponding Ramsey number $R(5,5)$.}\label{fig:R_5_5_40}
\end{figure}

We observe that reproducing lower bounds for Ramsey numbers with up to 28 vertices was possible. The first deviation occurs for the three-colour Ramsey number $R(3,3,4)$, where the algorithm was not able to find a clique-free configuration with 28 vertices (thus implying $R(3,3,4)\ge 29$ in Table~\ref{tab:multicolor}), but not 29 vertices (even though such a configuration is known to exist). For larger numbers of vertices, the results are less sharp. This is due to a significant increase in the required computation times when approaching the literature values. For instance, for $R(5,5)$, our algorithm took approximately 100s to converge to a clique-free 35-vertex configuration, about 1000s for 39 vertices, but growing to 40,000s for 40 vertices. The obtained clique-free configuration in the latter setup is shown in Figure~\ref{fig:R_5_5_40}. A similar steep incline of the required computation times occurred for multi-colour Ramsey numbers. In the outlook we discuss possible improvements that might overcome this impasse, since ultimately the goal is to match or possibly improve upon the known lower bounds.

%%%%%%%%%%%%%%%%%%%%%%%%%%%%%%%%%%%%%%%%%%%%%%%%%%%%%%%%%%%%%%%
\section{Challenges to the Monte Carlo approach}\label{sec:challange}
%%%%%%%%%%%%%%%%%%%%%%%%%%%%%%%%%%%%%%%%%%%%%%%%%%%%%%%%%%%%%%%
For several of the known optimal clique-free configurations the colourings are very peculiar. They are not obtained with a noisy computational approach, but by a more precise construction. As an example, we look at the known clique-free graph for 29 vertices~\cite{Kalbfleisch66,PiwakowskiRadziszowski98} for $R(3,3,4)$. If the vertices are numbered as $1,\ldots,29$, the colouring of the edges only depends on the ``distance" between the vertices, i.e.~the difference of the vertex labels (modulo the number of vertices). In other words, $A_{1,2}=A_{10,11}$, $A_{5,14}=A_{9,18}$ etc. This is called a cyclic colouring, for which the adjacency matrix takes a diagonal structure:
\begin{equation}
	\includegraphics[valign=c]{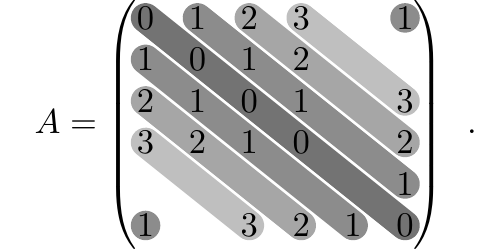}
\end{equation}
This non-local structure is unlikely to emerge from a Monte Carlo algorithm, which only updates the colourings of individual edges in each step. This may be why, when the number of vertices approaches the Ramsey number (or range) from below, the Monte Carlo approach slows down tremendously. Of course, it might be that clique-free configurations without such a non-local structure exist, but we did not manage to find such configurations for larger graphs, e.g.~the largest graph we obtained for  $R(3,3,4)$ had $N=28$ vertices.

\begin{figure}[t]
	\subfloat{\includegraphics[scale=1]{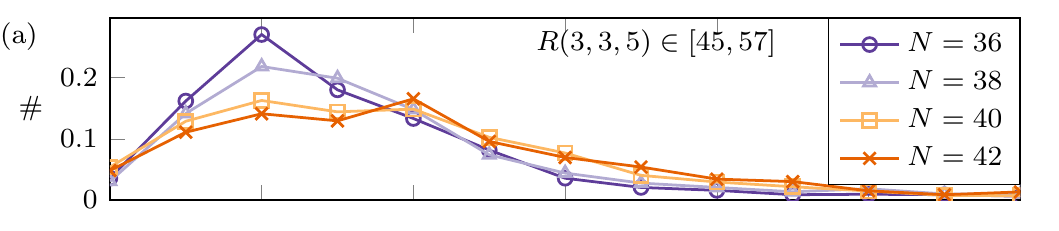}}\\
	\vspace{-0.5cm}
	\subfloat{\includegraphics[scale=1]{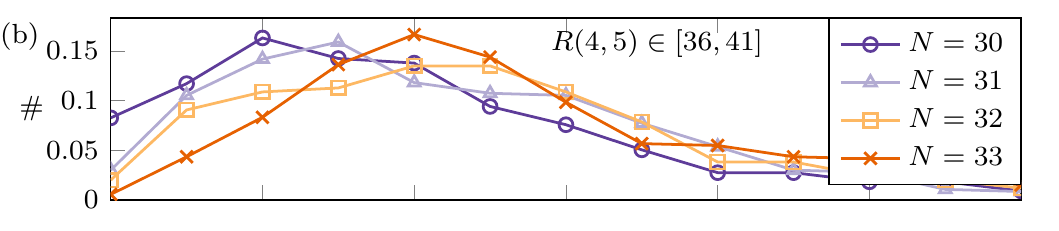}}\\
	\vspace{-0.7cm}
	\subfloat{\includegraphics[scale=1]{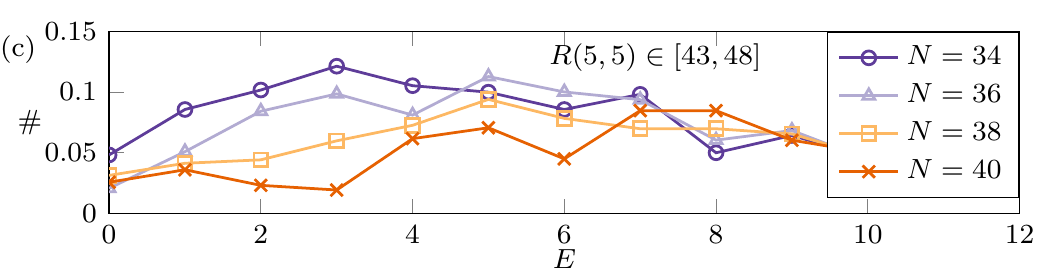}}\\	
	\caption{Normalised density of states connected to given clique-free configurations by changing the colour of individual edges. When increasing the number of vertices $N$ we observe a suppression at low energies, implying that fewer low-energy configurations can be found next to a clique-free configuration. This suppression is more pronounced for the two-colour Ramsey number.} \label{fig:MC_excitation_spectrum}
\end{figure}
To elaborate on the numerical challenge, we studied a physics-inspired quantity, namely a ``density of states". Specifically, for a given Ramsey number we considered several clique-free configurations for several graph sizes $N$. For each of them, we changed the colour on each edge (i.e.~we performed an inverse Monte Carlo step). We calculated the energy~(\ref{eq:energy}) of the resulting configuration and recorded their number as a function of their energy. The so-obtained number of configurations was then normalised to the total number of configurations connected to the initial clique-free one by changing the colour of a single edge, thus obtaining the density of states in the vicinity of clique-free configurations. 

The empirical results for three different Ramsey numbers are shown in Figure~\ref{fig:MC_excitation_spectrum}, where we have chosen $K_i=1$ for simplicity. We observe that when increasing the number of vertices, the density of states at low energy is suppressed. In other words, a typical configuration connected to the initial zero-energy state by a single colour flip will, on average, possess a higher energy. This implies that, in the inverse process, the final Monte Carlo step to find a clique-free configuration becomes more difficult, as fewer zero-energy states nearby exist. This is in line with our observation that the computation time increases as the number of vertices approaches the Ramsey number. There could be, of course, other factors contributing to the slowdown of the Monte Carlo search, like the decreasing overall number of clique-free configurations. 

%%%%%%%%%%%%%%%%%%%%%%%%%%%%%%%%%%%%%%%%%%%%%%%%%%%%%%%%%%%%%%%
\section{Conclusion and outlook}
%%%%%%%%%%%%%%%%%%%%%%%%%%%%%%%%%%%%%%%%%%%%%%%%%%%%%%%%%%%%%%%
In this article we have related the Ramsey problem to the task of finding zero-energy states of a suitably chosen statistical physics system. As a proof of principle, we treated the latter using well-known Monte Carlo methods, which gave us clique-free configurations on complete graphs, in turn establishing lower bounds on the corresponding Ramsey numbers. The results obtained were in a good agreement with literature values~\cite{Radziszowski21}, but we were unable to improve on any of the known lower bounds. 

We believe that there is much potential in modifications or fine tuning of the Monte Carlo algorithm, which could lead to novel insights into the Ramsey numbers. Several ideas come to mind: (i) We employed a Metropolis algorithm. Given the similarity of the problem to frustrated spin models and glassy systems, one may use alternative Monte Carlo approaches, such as parallel tempering~\cite{marinari_simulated_1992} or the Wang--Landau algorithm~\cite{dayal_performance_2004}. (ii) Given the peculiar form of the cyclic colourings, it might be better to perform the Monte Carlo updates on 'stripes' instead of individual edges. (iii) One could adapt the algorithm such that the number of $(x_i-1)$-cliques is maximised, in the hope that this opens 'room' for configurations without $x_i$-cliques. (iv) Avoiding the erratic energy landscape by adapting the energy functional~(\ref{eq:energy}) may speed-up the Monte Carlo search. 

Of course, one can try to minimise the energy functional using completely different methods, like the quantum evolution algorithm employed previously~\cite{GaitanClark12,Bian-13}. However, the intrinsic non-locality of the energy functional makes it impossible to employ standard one-dimensional techniques like density matrix renormalisation group method~\cite{White92,Schollwoeck11}.

%%%%%%%%%%%%%%%%%%%%%%%%%%%%%%
\section*{Acknowledgement}
%%%%%%%%%%%%%%%%%%%%%%%%%%%%%%
We would like to thank Stefan Wessel for useful discussions. This work is part of the D-ITP consortium, a program of the Dutch Research Council (NWO) that is funded by the Dutch Ministry of Education, Culture and Science (OCW). 
RJK is supported by a Vidi grant (639.032.614) of NWO.

\appendix
%%%%%%%%%%%%%%%%%%%%%%%%%%%%%%%%%%%%%%%%%%%%%%%%%%%%%%%%%%%%%%%
\section{Clique-biased Monte Carlo algorithm}\label{app:biasedbonds}
%%%%%%%%%%%%%%%%%%%%%%%%%%%%%%%%%%%%%%%%%%%%%%%%%%%%%%%%%%%%%%%
Recall that we are by no means interested in thermodynamic properties of the model defined by the energy functional~(\ref{eq:energy}). Instead we aim to find some zero-energy/clique-free state. We can therefore be biased with the choice of edges flipped in the Monte Carlo algorithm, potentially violating the condition of detailed balance. Only every $10\sim100$ steps we touch all the edges in the graph. In between the Monte Carlo algorithm only changes the colours of the edges in monochromatic cliques, since these are the sole changes that can lead a reduction in energy. Note that an $E>0$ state might not be related to any $E=0$ state by a simple swap of a edge from a monochromatic clique, instead it could be that the non-monochromatic part of the graph needs to be modified as well.

%%%%%%%%%%%%%%%%%%%%%%%%%%%%%%%%%%%%%%%%%%%%%%%%%%%%%%%%%%%%%%%
\section{Temperature variation}\label{app:temperature}
%%%%%%%%%%%%%%%%%%%%%%%%%%%%%%%%%%%%%%%%%%%%%%%%%%%%%%%%%%%%%%%
\begin{figure}[t]
	\centering
	\includegraphics[width=0.6\textwidth]{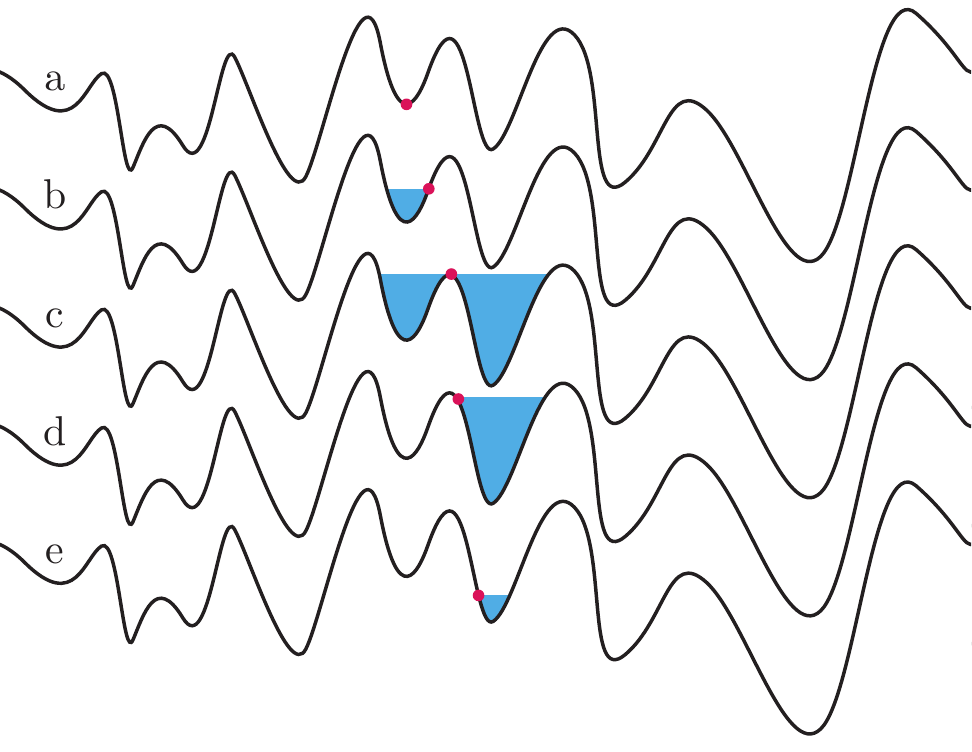}
	\caption{Schematic sketch of the procedure to vary the temperature in order to avoid sticking to local minima of the energy landscape (see text).}\label{fig:EL}
\end{figure}
The energy landscape of the energy functional~(\ref{eq:energy}) shows many local minima; see a sketch for illustration shown in Figure~\ref{fig:EL}. Thus simply lowering the temperature, as is done in a standard annealing procedure, will often cause the system to get stuck in one of these minima. To avoid this, we have slightly altered the approach. If the system gets stuck in a local minima [see panel (a)], i.e.~the algorithms stays at the same energy for many Monte Carlo steps, the temperature is slowly increased [see (b)]. This gives the system the chance to climb out of the local minimum [see (c)] and possibly descend into another minimum [see (d)]. As soon as the energy decreases again, the temperature is lowered [see (e)] to allow the system to descend to the bottom of the particular minimum. In this way the Monte Carlo algorithm can explore many local minima, ideally ending up at a zero-energy state. In practice, the temperature adaptation was a bit more elaborate, to let the system escape fractional minima (minima within minima\ldots) as well. Still, roughly the same ideas apply.

%%%%%%%%%%%%%%%%%%%%%%%%%%%%%%%%%%%%%%%%%%%%%%%%%%%%%%%%%%%%%%%
\section{Reducing for-loop calculation}\label{app:forloopreduction}
%%%%%%%%%%%%%%%%%%%%%%%%%%%%%%%%%%%%%%%%%%%%%%%%%%%%%%%%%%%%%%%
The number of fully connected cliques scales unfavourably with the number of vertices. As a result, explicitly checking all possibilities is computationally costly as it scales as $\sum_i {N\choose x_i}>\sum_i N^{x_i}/(x_i!)$. Fortunately we can, in the process, eliminate a lot of cliques and thus reduce the number to be checked. We implemented three reductions to speed up the code: 
\begin{enumerate}
	\item First, given $x_i$ we only consider cliques $Q_j^{x_i}=\{p_1,p_2,\ldots,p_{x_i}\}$ with $p_1<p_2<\ldots<p_{x_i}$, which covers all subgraphs and avoids overcounting due to the indistinguishability of the vertices.
	\item Second, we break for-loops early that are guaranteed to give no contribution to the energy. For a general $R(x_1,\ldots,x_l)$, the first two for-loops have given us the first colour $A_{p_1,p_2}$, fixing the colour of any monochromatic subgraph containing $p_1$ and $p_2$. The third for-loop yields $p_3$. If either $A_{p_1p_3}$ or $A_{p_2p_3}$ is not of the same colour, we can break the loop and continue with different $p_1$ and $p_2$. If $A_{p_1p_3}=A_{p_2p_3}=A_{p_1p_2}$, we continue with $p_4$, checking $A_{p_kp_4}$ for $k=1,2,3$, with the same conditional breaking. Note that, of course, we have to respect the subgraph sizes $x_i$ for each colour $i$.
	\item Third, in the Monte Carlo algorithm the colour-change on the edges is determined by the energy difference of the old and new state. Since the only energy difference can arise from the subgraphs that have the particular edge in common, we only have to check these subgraphs. This saves two for-loops.
\end{enumerate}

%%%%%%%%%%%%%%%%%%%%%%%%%%%%%%


\section*{References} 

\begin{thebibliography}{10}
\bibitem{Ramsey30}
F. P. Ramsey, \emph{On a problem of formal logic}, Proceedings of the London Mathematical Society \textbf{30}, 264 (1930); https://doi.org/10.1112/plms/s2-30.1.264.

\bibitem{graham1991}
R. L. Graham, B. L. Rothschild and J. H. Spencer, \emph{Ramsey Theory} (Wiley-Interscience Series in Discrete Mathematics and Optimization, 2nd Edition, 1991).

\bibitem{Radziszowski21}
S. Radziszowski, \emph{Small Ramsey numbers}, Electronic Journal of Combinatorics, Dynamic Surveys, DS1: Jan 15, 2021 (2021); https://doi.org/10.37236/21.

\bibitem{greenwood_combinatorial_1955}
R.~E. Greenwood and A.~M. Gleason.
\newblock \emph{Combinatorial {Relations} and {Chromatic} {Graphs}},
\newblock {Canadian Journal of Mathematics} \textbf{7}, 1 (1955); https://doi.org/10.4153/CJM-1955-001-4.

\bibitem{NewmanBarkema99}
M. E. J. Newman and G. T. Barkema, \emph{Monte Carlo methods in statistical physics} (Oxford University Press, Oxford, 1999).

\bibitem{GaitanClark12}
F. Gaitan and L. Clark, \emph{Ramsey numbers and adiabatic quantum computing}, Phys. Rev. Lett. \textbf{108}, 010501 (2012); https://doi.org/10.1103/PhysRevLett.108.010501.

\bibitem{Exoo1989OnTC}
G. Exoo, \emph{On Two Classical Ramsey Numbers of the Form R(3, n)}, SIAM J. Discret. Math., 488-490, 1989; https://doi.org/10.1137/0402043.

\bibitem{Bian-13}
Z. Bian, F. Chudak, W. G. Macready, L. Clark and F. Gaitan, \emph{Experimental determination of Ramsey numbers}, Phys. Rev. Lett. \textbf{111}, 130505 (2013); https://doi.org/10.1103/PhysRevLett.111.130505.

\bibitem{haanpaa_computational_2007}
H. Haanpää,
\newblock \emph{Computational Methods for Ramsey Numbers}, (Helsinki University of Technology, Laboratory for Theoretical Computer Science,
\newblock Research report 65, 2000).

\bibitem{Spencer94}
J. Spencer, \emph{The probabilistic method}, in \emph{Ten lectures on the probabilistic method} (Society for Industrial and Applied Mathematics, 1994); https://doi.org/10.1137/1.9781611970074.ch1.

\bibitem{Kalbfleisch66}
J.-G. Kalbfleisch, \emph{Chromatic graphs and {Ramsey}'s theorem}, (Ph.D. thesis, University of Waterloo, 1966).

\bibitem{PiwakowskiRadziszowski98}
K. Piwakowski and S. Radziszowski, \emph{$30\le R(3,3,4) \le 31$}, Journal of Combinatorial Mathematics and Combinatorial Computing \textbf{27}, 135 (1998).

\bibitem{White92}
S. R. White, \emph{Density matrix formulation for quantum renormalization groups}, Phys. Rev. Lett. \textbf{69}, 2863 (1992); https://doi.org/10.1103/PhysRevLett.69.2863.

\bibitem{Schollwoeck11}
U. Schollw\"ock, \emph{The density-matrix renormalization group in the age of matrix product states}, Ann. Phys. \textbf{326}, 96 (2011); https://doi.org/10.1016/j.aop.2010.09.012.

\bibitem{marinari_simulated_1992}
E.~Marinari and G.~Parisi,
\newblock \emph{Simulated tempering: {A} {new} {Monte} {Carlo} {scheme}},
\newblock {Europhys. Lett. (EPL)} \textbf{19}, 451 (1992); https://doi.org/10.1209/0295-5075/19/6/002.

\bibitem{dayal_performance_2004}
P.~Dayal, S.~Trebst, S.~Wessel, D.~Würtz, M.~Troyer, S.~Sabhapandit and S.~N.
Coppersmith, 
\newblock \emph{Performance limitations of flat-histogram methods},
\newblock {Phys. Rev. Lett.} \textbf{92}, 097201 (2004); https://doi-org.proxy.library.uu.nl/10.1103/PhysRevLett.92.097201.

\end{thebibliography}
\end{document}